\documentclass{amsart}
\usepackage{amssymb}
\usepackage{graphicx}
\def\({\left(}
\def\){\right)}

\newtheorem{lema}{Lemma}[section]

\newtheorem*{teorema*}{Theorem}

\newtheorem{remark}[lema]{Remark}

\newtheorem{lemma}{Lemma}[section]

\newtheorem{theorem}[lema]{Theorem}

\newtheorem{proposition}[lema]{Proposition}

\newtheorem{definition}[lema]{Definition}

  {\par \hfill \fbox{}}
  {\par \hfill \fbox{}}
  {\par \hfill }
  {\par \hfill }

\def\beq{\begin{equation}}
\def\eeq{\end{equation}}
\def\beginpf{\noindent{\bf Proof.} \quad}
\def\endpf{\qquad\hfill\rule{2.2mm}{2.2mm}\break}

\def\epsilon{\varepsilon}

\begin{document}

\title[Quaternionic Fock space]{Quaternionic Fock space  in an Euclidean unit ball $\mathbb{B}^n$}
\author{Sanjay Kumar}
\address{Department of Mathematics, Central University of Jammu,
Jammu 180 011, INDIA.} \email{sanjaykmath@gmail.com}
\author{Khalid Manzoor}
\address{Department of Mathematics, Central University of Jammu,
Jammu 180 011, INDIA.}
\email{khalidcuj14@gmail.com}


\date{\today.}

\begin{abstract}
In this paper, we  study the quaternionic counterpart of complex Fock spaces $\mathfrak{F}_{\alpha}^p  ( 0<p<\infty$  and for some parameter $\alpha$)  of entire  slice hyperholomorphic functions in an Euclidean unit ball $\mathbb{B}^n$ in $\mathbb{H}^n.$ We also study some basic properties of these  spaces.
\end{abstract}

\maketitle

\section{Introduction}
Consider  the 4-dimensional real algebra $\mathbb{H}=\{x_0+x_1i+x_2j+x_3k: x_l \in \mathbb{R}~\mbox{for} ~0\leq l\leq 3\}$  of  quaternions.
Then  $\{1, i, j, k\}$ form the standard basis of  quaternions with   imaginary units,  where   
 $Re(q)=x_0,\;\; Im(q)=x_1i+x_2j+ x_3k$  denote the  real part and  imaginary part   of $q$  respectively.
The conjugate  of $q$  is given by ~ $\bar q=Re(q)-Im(q)=x_0-(x_1i+x_2j+x_3k).$
  By   $|q| =\sqrt{q\bar q}=\sqrt{\sum_{l=0}^3}x_l^{2},~ \mbox{for}~ x_l \in \mathbb{R},$  we mean  the Euclidean norm on  $ \mathbb{H}.$  For any   non-zero quarternion $q$, its   multiplicative inverse  is given by $\displaystyle\frac{\bar q}{|q|^2}.$  The set $\mathbb{S}=\{q\in \mathbb{H}: q=x_1i+x_2j+x_3k  ~\mbox{and}~ x_1{^2}+x_2{^2}+x_3{^2}=1\}$  denotes  the two-dimensional  unit sphere of purely imaginary quaternions.  Any element $I\in \mathbb{S}$ is such that  $I^2=-1.$ This implies that the elements of S are  imaginary units.  The quaternion is considered as the union of  complex  plane  $\mathbb{C}_I=\mathbb{R}+\mathbb{R}I$  (also called slices),  each one is identified by an imaginary unit $I\in \mathbb{S}.$ Let   $\Omega_I= \Omega \cap \mathbb{C}_I,$ for  some 
domain $\Omega$ of $\mathbb{H}.$
 For any  quaternion $q$ we can write $q=x_0+x_1i+x_2j+x_3k=x_0+Im(q)=x_0+|Im(q)|I_q=x+yI_q,$
 with
$ I_q=\displaystyle\frac{Im(q)}{|Im(q)|}$ if  $|Im(q)|\neq 0,$ otherwise we take arbitrary $I$ in $\mathbb{S}.$ 

  For some more details on  slice holomorphic functions one can refer to \cite{fcs13, colo15} and references therein. In recent times the notion  of slice holomorphic functions has been studied intensively. The Hardy spaces, Bloch spaces, Dirichlet spaces, etc.  are studied in the slice hyperholomorphic settings. For Fock spaces in the slice hyperholomorphic settings, see \cite{alpa14}.

For details about recent work on slice hyperholomorphic and their applications to Operator Theory, Schur Analysis, Quantum Physics, etc one can refer to an excellent recent survey \cite{adle95, alpa12,fcs13,colo15, part92} and references therein.
For study of Fock spaces of holomorphic functions, we can refer to \cite{zhu12}.

 The work of this  paper is motivated by   \cite{alpa14, marco} and recent work of Ueki \cite{Sei14, Sei16}.
Here, we begin with some  basic results  in the quaternionic-valued slice regular functions.

\begin{definition}
Let $\Omega$ be a domain in $\mathbb{H}$.  A real differentiable function $f:\Omega \to \mathbb{H}$ is said to be the (left) slice regular or slice hyperholomorphic if for any $I\in \mathbb{S}, f_I$  is holomorphic in $\Omega_I$, i.e.,  $$ \left(\frac{\partial }{\partial x}+I\frac{\partial }{\partial y}\right)f_I(x+yI)=0,$$  where $f_I$ denote the restriction of $f$ to 
  $\Omega_I.$ The  class of slice regular functions on $\Omega$ is denoted  by $SR(\Omega).$
\end{definition}

For slice regular functions, we have the following useful result.
\begin{theorem} \cite[Theorem  2.7]{gent07}
A function $f: \mathbb{B} \to \mathbb{H}$ is said to be slice regular if and only if it has  a power series of the form 
\beq\label{eq:19}
f(q)=\displaystyle\sum_{n=0}^\infty q^na_n, ~~\mbox{where}~~a_n=\frac{1}{n!}\frac{\partial^n f(0)}{\partial x^n}
\eeq
 converging uniformly on $\mathbb{B}.$
\end{theorem}


\begin{lemma}\label{eq:1}
 \cite[Lemma 4.1.7 ]{fcs13}(Spliting Lemma) If $f$ is a slice regular function on the domain $\Omega,$ then for any 
$I,J\in  \mathbb{S},$  with $I\bot J$ there exist two holomorphic functions $F,L:\Omega_I =\Omega \cap \mathbb{C}_I\to\mathbb{C}_I $ 
such that 
\beq \label{eq:2}
f_I(z)=F(z)+L(z)J; ~~\mbox{for any }~~z=x+yI \in \Omega_I.
\eeq
\end{lemma}
\noindent One of the most important property of the slice regular functions is their Representation Formula.  It only  holds on the open sets which are stated below.
\begin{definition}
 Let $\Omega$ be an open set in  $\mathbb{H}$. We say  $\Omega$ is axially symmetric if for any  $q=x+yI_q \in \Omega,$ all the elements $x+yI$  are contained in $\Omega$, for all $I\in  \mathbb{S}$ and $\Omega$ is said to be slice domain if $\Omega\cap \mathbb{R}$ is non empty and  $\Omega\cap \mathbb{C}_I$ is a domain in $\mathbb{C}_I$ for all $I\in \mathbb{S}.$
\end{definition}
\begin{theorem} \label{eq:115}
  \cite[Theorem  4.3.2 ]{fcs13} (Representation Formula) Let  $f$ be a  slice regular function in the domain $\Omega\subset  \mathbb{H}.$ Then for any $J\in  \mathbb{S}$ and for all $z=x+yI \in \Omega$ such that $$f(x+yI)=\frac{1}{2}\{(1+IJ)f(x-yI)+(1-IJ)f(x+yI)\}.$$
\end{theorem}
\begin{remark}
Let $I,J$ be orthogonal imaginary units in $\mathbb{S}$ and $\Omega$ be  an  axially symmetric slice domain. Then the Splitting Lemma and the  Representation formula  generate a class of  operators on the slice regular functions as follows:

$$Q_I:SR(\Omega)\to hol({\Omega}_I)+hol({\Omega}_I)J$$ $$Q_I: f\mapsto  f_1+f_2J$$

$$P_I: hol({\Omega}_I)+hol({\Omega}_I)J \to SR(\Omega)$$
$$P_I[f](q)=P_I[f](x+yI_q)=\frac{1}{2}[(1-II_q)f(x+yI)+(1+II_q)f(x-yI)],$$ where   $hol({\Omega}_I)$ is the space of all holomorphic functions on  $\Omega_I$ for some $I \in \mathbb{S}.$ Also, $$P_I\circ Q_I=\emph{I}_{SR(\Omega)}~ \mbox{and}~Q_I\circ P_I=\emph{I}_{SR( hol({\Omega}_I)+hol({\Omega}_I))},$$ where $\emph{I}$ is an identity operator.

\end{remark}

Since pointwise product of functions does not preserve slice regularity, a new multiplication operation for regular functions is defined. In the special case of power series, the regular product (or $\star-$product) of $f(q) = \sum_{n= 0}^{\infty}q^{n}a_{n}$ and  $g(q) = \sum_{n= 0}^{\infty}q^{n}b_{n}$ is $$ f \star g(q) =  \sum_{n \ge 0} q^{n} \sum_{k = 0}^{n} a_{k}b_{n-k}. $$
The $\star-$product is related to the standard pointwise product by the following formula.
\begin{theorem} 
  \cite[Proposition   2.4 ]{arco14}
Let $f, g$ be regular functions on $\mathbb{B}.$ Then
$f \star g(q) = 0 $ if $f(q)  = 0 $ and $ f(q) g(f(q)^{-1} qf(q))$ if $ f(q) \neq 0.$ 
The reciprocal $f^{-\star}$ of a regular function $f (q) =\sum_{n= 0}^{\infty}q^{n}a_{n}$
with respect to the $\star-$product is
$$f^{\star}(q) = \frac{1}{f \star f^{c}(q)} f^{c}(q), $$ where $f^{c}(q) =\sum_{n= 0}^{\infty}q^{n}\overline{a_{n}} $ is the regular conjugate of $f.$ The function $f^{-\star}$
is regular on $\mathbb{B}\setminus(q \in \mathbb{B} | f \star f^{c} (q) = 0) $ and $f  \star f^{-\star} = 1$ there.
\end{theorem}

\section{Fock space}
\label{sec:2}
 Let $\mathbb{C}^n$ be  n-dimensional Euclidean  space and $dA$ be the normalized area measure on  $\mathbb{C}^n.$
 For $\alpha>0$ and $0<p< \infty,$ a holomorphic  function $f: \mathbb{C}^n \to  \mathbb{C}$ is said to be in Fock space  $\mathfrak{F}_{p, \mathbb{C}}$ if $$\displaystyle\left(\frac{\alpha}{\pi}\right)^n\int_{\mathbb{C}^n}\left|f(z)e^{\frac{-\alpha}{2}|z|^2}\right|^p dA(z) <\infty,$$
where  $ z=(z_1, z_2......,z_n)$ in  $\mathbb{C}^n,$ $dA_n(z) =(\frac{1}{\pi})^n\prod_{k=1}^n dx_kdy_k,~z_k=x_k+jy_k,~x_k, y_k \in \mathbb{R}.$ For 	$z=(z_1, z_2,....,z_n)$ and $w=(w_1, w_2,....,w_n)$ in $\mathbb{C}^n,$ we write $\langle z,w\rangle_\alpha$ for the Euclidean inner product $\displaystyle\sum_{k=1}^n z_k \overline{ w_k}$ and $|z|=\sqrt{\langle z,z\rangle}.$
 Let $\mathbb{B}^n(0,1)=\mathbb{B}^n=\{q=x+yI_q=(q_1, q_2,......,q_n): \displaystyle\sum_{k=1}^{n} |q_k|^2<1\}$ be the quaternionic unit ball  centered at origin in $\mathbb{H}^n$ and $\mathbb{B}^n\cap \mathbb{C}^n_I=\mathbb{B}_I^n$ denote n-dimensional  unit disk in the complex plane   $\mathbb{C}_I^n$ for $ I\in \mathbb{S}.$   
 Now we begin with  the following definitions.
\begin{definition}
 For $0<p<\infty,\alpha >0$ and some  $I\in \mathbb{S},$  the quaternionic right linear space of  entire slice regular functions $f$ is said to be the quaternionic slice  regular Fock  space on the unit ball $\mathbb{B}^n,$ if for any $q=x+yI_q\in\mathbb{B}^n $
$$\displaystyle\left(\frac{\alpha}{\pi}\right)^n\sup_{I\in \mathbb{S}}\int_{\mathbb{B}_I^n}\left|f(q)e^{\frac{-\alpha}{2}|q|^2}\right|^p dA_I(q) <\infty,$$that is,
$$\mathfrak{F}_{\alpha}^p=\left\{f\in SR( \mathbb{B}^n,\mathbb{H}^n):\left(\frac{\alpha}{\pi}\right)^n\sup_{I\in \mathbb{S}}\int_{\mathbb{B}_I^n}\left|f(q)e^{\frac{-\alpha}{2}|q|^2}\right|^p dA_I(q) \right\}<\infty,$$ where $dA_I(q)=\frac{1}{\pi}dx dy$
 and is M\"{o}bius invariant measure on $\mathbb{B}^n$.  If we write  $$d\lambda_{\alpha, I}(q)=\displaystyle\left(\frac{\alpha}{\pi}\right)^ne^{-\alpha |q|^2}dA_I(q),$$ then $\mathfrak{F}_{\alpha}^2$ is a Hilbert space (complete inner product space) with inner product  $\langle . , . \rangle_{\alpha}$ defined as $$\langle f, g\rangle_{\alpha}=\displaystyle \int_{\mathbb{B}_I^n}f(q)  \overline {g(q)}d\lambda_{\alpha, I}(q) $$with  the norm given  by
  $$\|f\|_{\mathfrak{F}_{\alpha}^p}= \left(\frac{\alpha}{\pi}\right)^{np}\sup_{I\in \mathbb{S}}\left(\int_{\mathbb{B}_I^n}\left|f(q)e^{\frac{-\alpha}{2}|q|^2}\right|^p dA_I(q)\right)^{\frac{1}{p}}.$$
\end{definition}
 By $\mathfrak{F}_{\alpha, I}^p,$ we denote the quaternionic right linear space of entire slice regular functions on  $\mathbb{B}^n$ such that 
 $$\displaystyle\left(\frac{\alpha}{\pi}\right)^n\int_{\mathbb{B}_I^n}\left|f(z)e^{\frac{-\alpha}{2}|z|^2}\right|^p dA_I(z) <\infty.$$ Furthermore, for each function $f\in \mathfrak{F}_{\alpha,I}^p,$ 
we define 
$$\|f\|_{\mathfrak{F}_{\alpha,I}^p} = \displaystyle\left(\frac{\alpha}{\pi}\right)^{np}\left(\int_{\mathbb{B}_I^n} \left|f(z)e^{\frac{-\alpha}{2}|z|^2}\right|^pdA_I(z):z=x+yI\in  \mathbb{B}^n\cap\mathbb{C}^n_I\right)^{\frac{1}{p}}.$$


\begin{remark}\label{eq:100} \cite[P. 499 ]{marco}\label{eq:100}
Let $I\in \mathbb{S}$ be such that $J\bot I.$ Then there exist holomorphic functions  $f_1, f_2: \mathbb{B}_I^n \to \mathbb{C}_I$ such that $Q_I[f]=f_1+f_2J$ for some holomorphic map $Q_I[f]$ in complex variable $z=(z_1, z_2,....z_n) \in \mathbb{B}_I^n$.  Then, we have 
$$\begin{array}{ccl}
\left|f_l(z)e^{\frac{-\alpha}{2}|z|^2}\right|^p&\leq &\left|f(z)e^{\frac{-\alpha}{2}|z|^2}\right|^p \\
&\leq& 2^{max\{0, p-1\}}\left|f_1(z)e^{\frac{-\alpha}{2}|z|^2}\right|^p+2^{max\{0, p-1\}}\left|f_2(z)e^{\frac{-\alpha}{2}|z|^2}\right|^p.
\end{array}$$
The condition  $f\in \mathfrak{F}_{\alpha,I}^p$  is equivalent to $f_1$ and $f_2$ belonging to n-dimensional   complex Fock space.  
\end{remark}
We can easily prove the following result.
\begin{proposition}\label{eq:20}
Let  $I\in  \mathbb{S}$ and  $\alpha>0.$ Then $f\in  \mathfrak{F}_{\alpha,I}^p, p>1$ if and only if $f\in  \mathfrak{F}_{\alpha}^p.$ Moreover, the spaces $( \mathfrak{F}_{\alpha,I}^p,\|.\|_{ \mathfrak{F}_{\alpha,I}^p})$ and $( \mathfrak{F}_{\alpha}^p,\|.\|_{ \mathfrak{F}_{\alpha}^p})$ have equivalent norms. More precisely, one has  $$\|f\|_{\mathfrak{F}_{\alpha,I}^p}^p\leq \|f\|_{\mathfrak{F}_{\alpha}^p}^p\leq 2^p\|f\|_{\mathfrak{F}_{\alpha,I}^p}^p.$$
\end{proposition}
\begin{proposition}
Let  $I,J\in \mathbb{S}$ and  let $p>1,\alpha >0$ and $f\in SR( \mathbb{B}^n,\mathbb{H}^n).$ Then  $f\in \mathfrak{F}_{\alpha, I}^p$  if and only if   $f\in \mathfrak{F}_{\alpha, J}^p.$
\end{proposition}
\beginpf
Since  $f\in SR(\mathbb{B}^n,\mathbb{H}^n)$ and $w=x+yJ \in \mathbb{B}^n_J$ and $z=x+yI \in \mathbb{B}^n_I,$ with $|w|=|z|.$  By Representation Formula, we have  $$f(w)=f(x+yJ)=\frac{1}{2}\displaystyle \left|(1-JI)f(z)+(1+JI)f(\bar z)\right|\leq |f(z)|+|f(\bar z)|.$$ This implies that
$$\begin{array}{ccl}
\displaystyle\int_{\mathbb{B}^n_J}\left|f(w)e^{\frac{-\alpha}{2}|w|^2}\right|^p dA_J(w)
&\leq&
2^{max\{p-1,0\}}\displaystyle\int_{\mathbb{B}^n_I}\left|f(z)e^{\frac{-\alpha}{2}|z|^2}\right|^p dA_I(z)\\
&+&\displaystyle 2^{max\{p-1,0\}}\int_{\mathbb{B}^n_I}\left|f(\bar z)e^{\frac{-\alpha}{2}|\bar z|^2}\right|^p dA_I(\bar z).\\
\end{array}$$ As $\bar z \to z,$ we have 
$$\begin{array}{ccl}
\displaystyle\left(\frac{\alpha}{\pi}\right)^n\displaystyle\int_{\mathbb{B}^n_J}\left|f(w)e^{\frac{-\alpha}{2}|w|^2}\right|^p dA_J(w)&\leq&2^{max\{p,1\}}\left(\displaystyle\left(\frac{\alpha}{\pi}\right)^n\displaystyle\int_{\mathbb{B}^n_I}\left|f(z)e^{\frac{-\alpha}{2}|z|^2}\right|^p dA_I(z)\right).
\end{array}$$ Thus, for  any $f\in \mathfrak{F}_{\alpha, I}^p,$ we see $f\in \mathfrak{F}_{\alpha, J}^p.$ By interchanging the roles of $I$ and $J,$  we obtain other result.
\endpf 

\begin{proposition}\label{eq:153}
 Suppose $p>1$ and $\alpha >0.$ If  $f$ is   in $ SR(\mathbb{B}^n,\mathbb{H}^n),$ then following assertations are equivalent:
\begin{itemize}
\item[(a)] $f\in \mathfrak{F}_{\alpha}^p(\mathbb{B}^n);$
 \item[(b)] $f\in \mathfrak{F}_{\alpha,I}^p(\mathbb{B}^n_I)$
 for some $I\in \mathbb{S}.$
\end{itemize}
\end{proposition}
\beginpf
To prove this, it is sufficient to  show  $(b)\Rightarrow (a).$  Suppose  $f\in \mathfrak{F}_{\alpha,I}^p.$ Let $q=x+yI_q$ and  $z=x+yI.$ Applying  Representation Formula  followed by triangle inequality and the fact that $|q|=|z|=|\bar z|,$ we see that
$$f(q)=\frac{1}{2}\left|(1-I_qI)f(x+yI)+(1+I_qI)f(x-yI)\right|\leq |f(z)|+|f(\bar z)|.$$ This implies that
$$\begin{array}{ccl}
\displaystyle\left(\frac{\alpha}{\pi}\right)^n\displaystyle\int_{\mathbb{B}^n_I}\left|f(q)e^{\frac{-\alpha}{2}|q|^2}\right|^p dA_I(q)
&\leq&2^{max\{p-1,0\}}\displaystyle\left(\frac{\alpha}{\pi}\right)^n\displaystyle\int_{\mathbb{B}^n_I}\left|f(z)e^{\frac{-\alpha}{2}|z|^2}\right|^p dA_I(z)\\
&+&\displaystyle 2^{max\{p-1,0\}} \left(\frac{\alpha}{\pi}\right)^n\int_{\mathbb{B}^n_I}\left|f(\bar z)e^{\frac{-\alpha}{2}|\bar z|^2}\right|^p dA_I(\bar z)\\
&\leq&2^{max\{p,1\}}\displaystyle\left(\frac{\alpha}{\pi}\right)^n\displaystyle\int_{\mathbb{B}^n_I}\left|f(z)e^{\frac{-\alpha}{2}|z|^2}\right|^p dA_I(z)\\
&<&  \infty.
\end{array}$$ 
Thus, the  condition $(a)$ holds.
\endpf

\begin{remark} By $L^{p}(\mathbb{B}_I, d\lambda_{\alpha, I}, \mathbb{H})$  we define   the set of functions $g: \mathbb{B}^n_I \to \mathbb{H} $ such that $$\displaystyle\int_{\mathbb{B}^n_I} |g(w)|^p d\lambda_{\alpha, I}(w)<\infty,$$ where $d\lambda_{\alpha, I}(w)=\frac{\alpha}{\pi}e^{-\alpha|z|^2}dA_I(w),$ for $\alpha>0$  is called  the Gaussian probability measure. Note that for $J\in \mathbb{S}$ with $J\bot I$ and $g=g_1+g_2J$ with $g_1,g_2:\mathbb{B}^n_I\to \mathbb{C}_I,$ then $g\in L^{p}(\mathbb{B}_I, d\lambda_{\alpha, I}, \mathbb{H})$ if and only if  $g_1 ,g_2 \in L^{p}(\mathbb{B}_I, d\lambda_{\alpha, I}, \mathbb{C}_I).$\\ Clearly, $\mathfrak{F}_{\alpha}^p$ is closed subspace of  $L^p(\mathbb{B}_I, d\lambda_{\alpha, I}, \mathbb{H}).$
In complex analysis, the reproducing kernel of complex Fock space for $p=2$ is given by
 $$K_{\alpha}^{\mathbb{C}_I^n}(z,w)=e^{\alpha\langle z,w\rangle};~~z,w\in \mathbb{C}^n_I.$$
\end{remark}
This gives the motivation  for the following definition.

\begin{definition}
For any $q\in \mathbb{B}^n,$ the slice regular exponential function is given by $$e^q= \sum_{n=0}^{\infty} \frac{q^n}{n!}.$$
Let  $\displaystyle e^{zw}= \sum_{n=0}^{\infty} \frac{z^n w^n}{n}$ be a holomorphic function in variable $z$ in the complex plane $\mathbb{C}^n _I$. Clearly, $e^{zw}$ is not slice regular in both  variable.  Setting $\displaystyle e^{qw}_{\star}= \sum_{n=0}^{\infty} \frac{q^n w^n}{n!},$ then we see that  the function  is left slice regular in q and right slice regular in w, where $\star$  denote  the slice regular  product. By Representation formula, we can obtain the extension  of  function $e^{zw}$ to $\mathbb{H},$  as $$ext(e^{zw})=\frac{1}{2}\{(1-IJ)e^{zw}+(1+IJ)e^{\bar{z}w}\}=e^{qw},$$ where $q\in \mathbb{B}^n$ and for some arbitrary $w.$
For $I_q=\frac{Im(q)}{|Im(q)|}\in \mathbb{S}$ and $\alpha>0,$ we define $$B_{\alpha}(q,w)=e^{\alpha q\bar w}_{\star}~~ for~ each ~q=x+yI_q\in \mathbb{B}^n,$$ and  is called slice regular reproducing kernel  of quaternionic  Fock space. 
\end{definition}

\noindent The following theorem is a  quaternionic version of   \cite[Theorem 3.1]{wall89}. 
\begin{theorem}\label{eq:103}
Let $I\in \mathbb{S}.$ Then for  any $p>1,$ there exists a sequence  $\{z_k\}_{k>1}$ in $\mathbb{B}^n_I$ with the following property:
An entire function $f\in  \mathfrak{F}_{\alpha}^p$ if and only if  
\beq\label{eq:101}
f(z)=P_{I}\sum_{k=1}^\infty   w_{z_k}(z)a_k ,
\eeq
 where $\{a_k\}\in l^p(\mathbb{H}^n)$ and  $w_{z_k}(z)=e_{\star}^{\alpha z \bar {z_k}-\frac{\alpha}{2}|z_k|^2}$  is the slice regular  normalized reproducing kernel with  $P_{I}=\frac{\vec z}{\vec {\|z\|}}.$
\end{theorem}
\beginpf
 Let $J\in \mathbb{S}$ be such that $J\bot I.$  Let $f\in \mathfrak{F}_{\alpha}^p.$  Then $f\in \mathfrak{F}_{\alpha,I}^{p}$ for $I\in \mathbb{S}.$  If we restrict $f$  on $\mathbb{B}^n_I,$  then  we have  $ f_I(z)=Q_I[f](z)=f_{1}(z)+f_{2}(z)J$ with  holomorphic functions $f_1, f_2 : \mathbb{B}^n \cap  \mathbb{C}^n_I\to \mathbb{C}_I.$  Applying \cite[Theorem 3.1]{wall89} to $f_1$ and $f_2,$  we  can find the   sequences $\{a_{1,k}\}_{k\geq 1}$ and $\{a_{2,k}\}_{k\geq 1}$  in $ l^p(\mathbb{C}_I^n)$  such that  $$ f_1(z)=\sum_{k=1}^\infty  w_{z_k}(z)a_{1,k}  ~~\mbox{and}~~ f_2(z)=\sum_{k=1}^\infty w_{z_k}(z)a_{2,k}.$$  So, we can write $$Q_I[f](z)=f_{1}(z)+f_{2}(z)J=\sum_{k=1}^\infty  w_{z_k}(z)a_{1,k} +\sum_{k=1}^\infty  w_{z_k}(z)a_{2,k}J.$$ As $P_I\circ Q_I[f]=$I$_{SR(\mathbb{B}_I^n)}$ an identity operator, we have  
$$\begin{array}{ccl}
f=P_I\circ Q_I[f](z)&=&P_I\displaystyle\left\{\sum_{k=1}^\infty  w_{z_k}(z) a_{1,k}+\sum_{k=1}^\infty   w_{z_k}(z)a_{2,k} J\right\}\\
&=&P_I \displaystyle\left\{\sum_{k=1}^\infty  w_{z_k}(z)\left(a_{1,k} +a_{2,k} J\right) \right\}\\
&=&P_I\displaystyle\left\{\sum_{k=1}^\infty   w_{z_k}(z)a_k\right\},
\end{array}$$
where $a_k=a_{1,k} +a_{2,k} J \in l^p(\mathbb{H}^n)$. Conversely, suppose  the condition  (\ref{eq:101}) holds. We claim $f\in  \mathfrak{F}_{\alpha}^p.$ Then for any $J \in \mathbb{S}, $  we have  $a_k=a_{1,k} +a_{2,k} J,$ where  $\{a_{l,k}\}_{k\geq 1} \in  l^p(\mathbb{C}^n_I),~l=1,2.$  So, we can write  $$Q_I[f](z)=\sum_{k=1}^\infty a_{1,k} w_{z_k}(z)+\sum_{k=1}^\infty a_{2,k} w_{z_k}(z)J.$$  Since  $\{a_k\}_{k\geq 1}$ lie in  $l^p(\mathbb{H}^n),$ it follows that  $\{a_{l,k}\}_{k\geq 1}$ belong to $l^p(\mathbb{C}^n_I).$ This implies that $f_1, f_2$ lie in the complex Fock space  $ \mathfrak{F}_{\alpha,I,\mathbb{C}_I}^{p}$  which is equivalent to  $f \in \mathfrak{F}_{\alpha,I}^{p}$ and hence $f\in \mathfrak{F}_{\alpha}^{p}.$  
\endpf\\
\noindent The  proof of the  following proposition is analogus to  \cite[Theorem 2.19 ]{marco}.\\

\begin{theorem}
Suppose $1\leq p<\infty$ and $\alpha >0.$ If   $f\in \mathfrak{F}_{\alpha}^p$ with some  $I\in \mathbb{S},$ then for any $z\in  \mathbb{B}^n_I,$ there exists some $M>0$ such that $$|f(z)|\leq 2^{max{\{p,1\}}} M\|f\|_{\mathfrak{F}_{\alpha}^p}.$$
\end{theorem}

\begin{proposition}
  Let $1< p<\infty$ and $\alpha >0.$ If   $\{f_n\}_{n\in \mathbb{N}}$ is any sequence  in $ \mathfrak{F}_{\alpha}^p(\mathbb{B}^n),$ then  $\{f_n\}$ converges weakly  to $f \in \mathfrak{F}_{\alpha}^p(\mathbb{B}^n)$ if and only if there exists $\lambda>0$ such that  $\|f_n\|_{\mathfrak{F}_{\alpha}^p}\leq \lambda$ and  $f_n \to f$ uniformly on the compacts subsets of $\mathbb{B}^n.$ 
\end{proposition}
\beginpf
Suppose    sequence $f_n$ is weakly convergent  in $\mathfrak{F}_{\alpha}^p$ and so in   $\mathfrak{F}_{\alpha,I}^p.$ Let $J\in \mathbb{S}$ be such that $J\bot I$ and let $f_{n,1},f_{n,2}$ be holomorphic functions such that $Q_I[f_n]=f_{n,1}+f_{n,2}J.$ By Remark \ref{eq:100},  $f_{n,1},f_{n,2}$ lie in the complex Fock space. Then, from corresponding result in the complex case, we see that   the sequences $\{f_{n,1}\}$  and  $\{f_{n,2}\}$   converges weakly to functions $f_1$ and $f_2$ as $n\to \infty,$ respectively if and only if  there exists some positive constant $M$ such that $\|f_{n,l}\|_{\mathfrak{F}_{\alpha, \mathbb{C}_I}^p }\leq M$ and $f_{n,l} \to f_l$ uniformly on the compact subsets of $\mathbb{B}^n_I,$ for $l=1,2$ if and only if   $\|f_{n}\|_{\mathfrak{F}_{\alpha,I}^p }\leq  2^{p-1}\|f_{n,l}\|_{\mathfrak{F}_{\alpha, \mathbb{C}_I}^p } \leq 2^{p-1}M=\lambda$ for some $\lambda=2^{p-1}M$ and for any $f=p_I[f_1+f_2J]$ in $\mathfrak{F}_{\alpha,I}^p,$ $$|f_n-f|\leq |f_{n,1}-f_1|+|f_{n,2}-f_2|\to 0, ~~as~~n\to \infty$$  uniformly on the compacts subsets of $\mathbb{B}^n_I.$   
\endpf
Now, we have the following definitions.
\begin{definition}
By space $\mathfrak{F}_{\infty,\alpha},$ we mean the space of all   slice regular Fock functions $f$ such that 
$$ \sup_{q\in\mathbb{B}^n} |f(q)|e^{\frac{-\alpha}{2}|q|^2}<\infty,$$
that is, $$\mathfrak{F}_{\infty,\alpha}=\displaystyle\left\{f\in  SR(\mathbb{B}^n,\mathbb{H}^n): \sup |f(q)|e^{\frac{-\alpha}{2}|q|^2}:~q=x+yI_q\in\mathbb{B}^n \right\}<\infty.$$ The space $\mathfrak{F}_{\infty,\alpha}$ is a Banach space under the norm defined as $$\|f\|_{\mathfrak{F}_{\infty,\alpha}}=\sup\left\{ |f(q)|e^{\frac{-\alpha}{2}|q|^2}:~q\in\mathbb{B}^n\right\}.$$
\end{definition}
\begin{definition}
We define the space  $\mathfrak{F}_{\infty,\alpha,I}$ as the quaternionic right linear space of all slice regular functions $f$ defined on 
$\mathbb{B}^n_I$ such that $$\|f\|_{\mathfrak{F}_{\infty,\alpha,I}}=\sup\left\{ |Q_I[f](z)|e^{\frac{-\alpha}{2}|z|^2}:~z=x+yI\in\mathbb{B}^n_I\right\}<\infty.$$
\end{definition}
\begin{definition}
The  space $\mathfrak{F}_{\infty,\alpha}^0$  is a closed subspace of $\mathfrak{F}_{\infty,\alpha}$  and  is defined as the space of all slice regular  functions in $\mathbb{B}^n$ such that  $$\displaystyle \lim_{|q|\to \infty}\left|f(q)\right|e^{\frac{-\alpha}{2}|q|^2}=0.$$
\end{definition}
{\begin{remark}
Let $I,J\in \mathbb{S}$ be such that $I\bot J$ and let $f\in SR(\mathbb{B}^n,\mathbb{H}^n).$ Furthermore, let $f_1, f_2:\mathbb{B}^n_I\to \mathbb{C}_I $ be holomorphic functions in the complex Fock space $\mathfrak{F}_{\infty,\alpha, \mathbb{C}_I}^{0}$ such that $Q_I[f]=f_1+f_2J.$ Then $$\displaystyle\left|f(z)e^{\frac{-\alpha}{2}|z|^2}\right|^2=\left|f_1(z)e^{\frac{-\alpha}{2}|z|^2}\right|^2+\left|f_2(z)e^{\frac{-\alpha}{2}|z|^2}\right|^2~\mbox{for all}~z\in \mathbb{B}^n_I.$$ By Representation Formula, it  follows that $f\in \mathfrak{F}_{\infty,\alpha}^0$ if and only if $f_1, f_2$ lie in  the complex space $\mathfrak{F}_{\infty,\alpha,\mathbb{C}_I}^{0}.$
\end{remark}
\noindent  By Representation Formula, we have the  following result.
\begin{proposition}\label{eq:21}
Suppose $I\in  \mathbb{S}$ and  $\alpha>0.$ Then $f\in  \mathfrak{F}_{\infty,\alpha,I}$ if and only if $f\in \mathfrak{F}_{\infty,\alpha}.$ Futhermore, the spaces $( \mathfrak{F}_{\infty,\alpha,I},~\|.\|_{\mathfrak{F}_{\infty,\alpha,I}})$ and $( \mathfrak{F}_{\infty,\alpha},~\|.\|_{ \mathfrak{F}_{\infty,\alpha}})$ have equivalent norms. More precisely, one has  $$\|f\|_{\mathfrak{F}_{\infty,\alpha,I}}\leq \|f\|_{\mathfrak{F}_{\infty,\alpha}}\leq 2\|f\|_{\mathfrak{F}_{\infty,\alpha,I}}.$$
\end{proposition}

\begin{proposition}
For $p>1$ and $\alpha>0$ with multi-index $m$ such that $m=(m_1, m_2,.......,m_n)$ and $|m|=|m_1|+ |m_2|+.......+|m_n|.$ Then for all $f\in  SR(\mathbb{B}^n,\mathbb{H}^n)$ and all $z\in \mathbb{B}^n_I,$ for some  $I\in \mathbb{S}$ $$\|a_mz^m\|_{\mathfrak{F}_{\infty,\alpha}}\leq 2^{max\{p,1\}} \prod_{k}\sqrt{\frac{m_k}{2}}\displaystyle\sup_{z\in  \mathbb{B}^n_I} |f(z)|e^{\frac{-\alpha}{2}|z|^2}.$$
\end{proposition}
\beginpf
Let  $J\in \mathbb{S}$ with $I\bot J.$ Then by Splitting Lemma, we can find  two  holomorphic functions   $f_1, f_2 \in  hol(\mathbb{B}_I^n)$  such that $Q_I[f]=f_1+f_2J.$  Therefore from \cite[Lemma 2]{Sei16}, we conclude that
$$\begin{array}{ccl}
\|a_mz^m\|_{\mathfrak{F}_{\infty,\alpha}}&\leq& 2^{max\{p-1,0\}}\displaystyle\left(\|a_{m,1}z^m\|_{\mathfrak{F}_{\infty,\alpha, \mathbb{C}_I}}+\|a_{m,2}z^m\|_{\mathfrak{F}_{\infty,\alpha,\mathbb{C}_I}}\right)\\
&\leq&2^{max\{p-1,0\}}\displaystyle\left(\prod_{k}\sqrt{\frac{m_k}{2}}\displaystyle\sup_{z\in  \mathbb{B}^n_I} |f_1(z)|e^{\frac{-\alpha}{2}|z|^2}+\prod_{k}\sqrt{\frac{m_k}{2}}\displaystyle\sup_{z\in  \mathbb{B}^n_I} |f_2(z)|e^{\frac{-\alpha}{2}|z|^2}\right)\\
&\leq&2^{max\{p-1,0\}}\displaystyle\prod_{k}\sqrt{\frac{m_k}{2}}\sup_{z\in  \mathbb{B}^n_I}(|f_1(z)|+|f_2(z)|)e^{\frac{-\alpha}{2}|z|^2}\\
&\leq& 2^{max\{p,1\}}\displaystyle\prod_{k}\sqrt{\frac{m_k}{2}}\sup_{z\in  \mathbb{B}^n_I}(|f(z)|)e^{\frac{-\alpha}{2}|z|^2}.\\
\end{array}$$ Hence the result.
\endpf
\begin{proposition}
 For  $0<r<1$, the space $ \mathfrak{F}_{\infty,\alpha}$ is the closure with respect to norm $ \|.\|_{ \mathfrak{F}_{\infty,\alpha}}$  of the set of quaternionic polynomials of the form 
$$P_k^{(r)}(q)=\sum_{l=1}^k \sum_{|m|=l}r^{|m|}q^ma_m~~for~all~~a_m\in \mathbb{H}^n.$$  In particular,  $ \mathfrak{F}_{\infty,\alpha}$ is separable.
\end{proposition}
\beginpf
 Let $ f\in \mathfrak{F}_{\infty,\alpha}.$  Then by  orthogonality of  imaginary units $I,J$ one can write $f_I=Q_I[f]=f_1+f_2J,$ where $f_1, f_2$ are holomorphic functions in the complex Fock space.  Since the polynomials are dense in the complex Fock space, from \cite[Lemma 3]{Sei16}, we find the complex polynomials of the form $$P_{k,1}^{(r)}(z)=\sum_{l=1}^k \sum_{|m|=l}r^{|m|}z^m\beta_{k,m}~~\mbox{and}~~P_{k,2}^{(r)}(z)=\sum_{l=1}^k \sum_{|m|=l}r^{|m|}z^m\gamma_{k,m}$$ such that $$\|f_n- P_{k,n}^{(r)}\|_{ \mathfrak{F}_{\infty,\alpha,\mathbb{C}_I^n}}\to 0,~\mbox{ as}~ k\to \infty~\mbox{for}~ n=1,2.$$ \\ Now, let $a_{k,m}=\beta_{k,m}+\gamma_{k,m}J.$ Then, we have $$P_k^{(r)}(q)=\sum_{l=1}^k \sum_{|m|=l}r^{|m|}q^m\left(\beta_{k,m}+\gamma_{k,m}J\right)=\sum_{l=1}^k \sum_{|m|=l}r^{|m|}q^ma_{k,m}.$$ Therefore,  by Remark \ref{eq:21}, it follows that
 $$\begin{array}{ccl}
\|f- P_{k}^{(r)}\|_{ \mathfrak{F}_{\infty,\alpha}}&\leq& 2\|f- P_{k}^{(r)}\|_{ \mathfrak{F}_{\infty,\alpha,I}}\\
&\leq& 2\|f_1- P_{k,1}^{(r)}\|_{ \mathfrak{F}_{\infty,\alpha, \mathbb{C}_I}}+ \|f_2- P_{k,2}^{(r)}\|_{ \mathfrak{F}_{\infty,\alpha, \mathbb{C}_I}}\\
&\to& 0,~as~k\to \infty.
\end{array}$$
Hence the space $\mathfrak{F}_{\infty,\alpha}$  is separable.
\endpf
\begin{proposition}
Suppose $f\in \mathfrak{F}_{\infty,\alpha},  \alpha>0.$ Then $f\in  \mathfrak{F}_{\infty,\alpha}^0$ if and only if 
\beq\label{eq:200}
\lim_{r\to 1}\|f_r-f\|_{ \mathfrak{F}_{\infty,\alpha}}=0,
\eeq
where the power series expansion of $f(z)=\sum q^ma_m$ and so we can write $$f_r(q)=f(rq)=\sum_{l=1}^k \sum_{|m|=l}r^{|m|}q^ma_m~for~all~a_m\in \mathbb{H}^n,~q\in \mathbb{B}^n~and~0<r<1$$ with multi-index $m=(m_1, m_2,.....,m_n)$ and $|m|=|m_1|+|m_2|+......+|m_n|.$
\end{proposition}
\beginpf
Let $I,J\in  \mathbb{S} $ be  such that $I\bot J.$ Then for any $a_m\in  \mathbb{H}^n,$ there exist $a_{m,1}$ and   $a_{m,2}$ in $ \mathbb{C}_I^n$ such that $a_m=a_{m,1}+a_{m,2}J.$ Now, $f\in \mathfrak{F}_{\infty,\alpha}$ implies $f\in \mathfrak{F}_{\infty,\alpha ,I}.$ So we can choose   holomorphic functions $f_1, f_2$ in the complex space $ \mathfrak{F}_{\infty,\alpha,\mathbb{C}_I}$  such that $Q_I[f]=f_1+f_2J.$ Furthermore, $f(z)=\sum z^m a_{m,1}+\sum z^ma_{m,2}J=f_1(z)+f_2(z)J.$ Thus, by Remark \ref{eq:21}  and for $n=1,2,$ we have
 $$\begin{array}{ccl}
\|f_{r,n}-f_n\|_{\mathfrak{F}_{\infty,\alpha, \mathbb{C}_I}}&\leq& \|f_{r}-f\|_{\mathfrak{F}_{\infty,\alpha,I}}\\
&\leq&  \|f_{r}-f\|_{\mathfrak{F}_{\infty,\alpha}}\\
&\leq&2\|f_{r}-f\|_{\mathfrak{F}_{\infty,\alpha,I}}\\
&\leq&2( \|f_{r,1}-f_1\|_{\mathfrak{F}_{\infty,\alpha, \mathbb{C}_I}}+\|f_{r,2}-f_2\|_{\mathfrak{F}_{\infty,\alpha, \mathbb{C}_I}}).
\end{array}$$
Thus because of  \cite[Theorem 2]{Sei16}, the condition  (\ref{eq:200}) holds if and only if $f_1, f_2 \in \mathfrak{F}_{\infty,\alpha,\mathbb{C}_I}^0$ which is equivalent to $f\in\mathfrak{F}_{\infty,\alpha}^0.$ 
\endpf

\noindent Here, we investigate the $t^{th}$-derivative criterion for the space $f\in\mathfrak{F}_{\infty,\alpha}^0$  on a unit ball  $\mathbb{B}$ in $ \mathbb{H}.$
\begin{definition}
For each $t\in \mathbb{N} \cup \{0\},$ we define the space $\mathfrak{F}^{(t)}$  of entire  slice regular  function $f$ such that $$\displaystyle \lim_{|q|\to \infty} \frac{\left|\displaystyle {\partial^{(t)}_{x_0} f}(q)\right|}{(1+|q|)^t}e^{\frac{-\alpha}{2}|q|^2}=0,~~q\in \mathbb{B}.$$

\end{definition}
\begin{remark}\label{eq:150}
\textnormal{Let $f\in SR(\mathbb{H},\mathbb{B})$  and  $I,J\in \mathbb{S}$ be  orthogonal imaginary units. Let $f_1, f_2$ be two holomorphic functions such that $Q_I[f]=f_1+f_2J.$ Then for all $z\in \mathbb{B}_I$ $$|{\partial^{(t)}_{x_0} f}(q)|\leq |f_1^{(t)}(z)|+|f_2^{(t)}(z)|.$$ By Theorem \ref{eq:115}, it follows that $f$ lie in  $\mathfrak{F}^{(t)}$  if and only if $f_l\in  \mathfrak{F}^{(t)}_{\mathbb{C}_I}$ for $l=1,2.$ }
\end{remark}
\noindent We can easily prove the following results.
\begin{theorem}
Suppose  $0\leq r<1,$ if  $f\in\mathfrak{F}^{(t)}, $ then 
$$\displaystyle \lim_{r\to 1}\|f_r-f\|_{ \mathfrak{F}_{\infty,\alpha}}=0,$$
 where
$f_r(q)=f(rq)=\displaystyle\sum_{k=0}^{\infty}r^kq^ka_k ~~for~all~~a_k\in \mathbb{H},~q\in \mathbb{B}.$
\end{theorem}
\begin{proposition}
$ \mathfrak{F}^{(t)}$ is the closure with respect to norm $ \|.\|_{ \mathfrak{F}_{\infty,\alpha}}$  of the set of quaternionic polynomials of the form 
$$P_k(q)=\sum_{m=0}^k r^m q^ma_m~~for~all~~a_m\in \mathbb{H},$$
for~ some ~positive ~integer~~k~and ~$0<r<1.$ In particular, $ \mathfrak{F}^{(t)}$ is separable.
\end{proposition}
\begin{theorem}
Let $f\in SR(\mathbb{B},\mathbb{H}).$ Then for any positive integer $t, f\in\mathfrak{F}_{\infty,\alpha}$ if and only if  
\beq\label{eq:151}
\frac{|{\partial^{(t)}_{x_0} f}(q)|}{(1+|q|)^t}e^{\frac{-\alpha}{2}|q|^2}\in l^{\infty}(\mathbb{H}).
\eeq
\end{theorem}
\beginpf
Let $f\in\mathfrak{F}_{\infty,\alpha}.$ Then  $f\in \mathfrak{F}_{\infty,\alpha,I}.$  Let $I,J\in $  with  $I\bot J.$ So, we can find holomorphic functions $f_1, f_2: \mathbb{B}\cap \mathbb{C}_I \to \mathbb{C}_I$ such that $f_I=f_1+f_2J.$ Now for any $z\in \mathbb{B}_I, $ we have  $$\frac{|f_l^{(t)}(z)|}{(1+|z|)^t}e^{\frac{-\alpha}{2}|z|^2}\leq \frac{|{\partial^{(t)}_{x_0} f}(z)|}{(1+|z|)^t}e^{\frac{-\alpha}{2}|z|^2}\leq \frac{|f_1^{(t)}(z)|}{(1+|z|)^t}e^{\frac{-\alpha}{2}|z|^2}+\frac{|f_2^{(t)}(z)|}{(1+|z|)^t}e^{\frac{-\alpha}{2}|z|^2}~\mbox{for}~l=1,2.$$ On applying  \cite[Theorem 1]{Sei14} to $f_1$ and $f_2,$ it follows that the condition (\ref{eq:151}) holds if and only if $ \displaystyle\frac{|f_l^{(t)}(z)|}{(1+|z|)^t}e^{\frac{-\alpha}{2}|z|^2}\in L^{\infty}(\mathbb{C}_I), ~l=1,2.$  Thus, we say that the above condition holds if and only if $f_1, f_2$ belong to complex space $\mathfrak{F}_{\infty,\alpha,\mathbb{C}_I}$ which is equivalent to $f\in \mathfrak{F}_{\infty,\alpha,I}$ and so $f\in \mathfrak{F}_{\infty,\alpha}.$
\endpf

\end{document}